\documentclass[12pt,a4paper,reqno]{amsart}
\usepackage[utf8]{inputenc}
\usepackage{amsmath}
\usepackage{amssymb}
\usepackage{tikz}
\usepackage{cases}
\usepackage{color}
\usepackage[left=3cm,right=2cm,top=3cm,bottom=2cm]{geometry}\usepackage[square,sort,comma,numbers]{natbib}

\newtheorem{Teo}{Theorem}

\newtheorem{Lem}[Teo]{Lemma}

\title[Complex invariant Einstein metrics]{Complex invariant Einstein metrics on $SO_{2(n_1+n_2+n_3)+1}/U_{n_1} \times U_{n_2} \times SO_{2n_3+1}$ \\ and Ricci-flat manifolds}

\author[A. Lavrov]{Aleksei Lavrov}

\begin{document}

\maketitle

\begin{abstract}
\noindent We prove that the number of complex invariant Einstein metrics on the flag manifold $M_{n_1,n_2,n_3}=SO_{2(n_1+n_2+n_3)+1}/U_{n_1} \times U_{n_2} \times SO_{2n_3+1}$ is equal to 132, except when the parameters $n_1, n_2, n_3$ satisfy one of some algebraic equations. Also the family of (real) non-flat Ricci-flat metrics on the Euclidean spaces will be constructed using the method of Inonu-Wigner contractions of Lie algebras.
\end{abstract}

\section{Introduction}

A Riemannian manifold $(M, g)$ is called Einstein if the Ricci tensor $\text{Ric}_{g}$ is a constant multiple of
the metric $g$, i.e. $\text{Ric}_{g} = \lambda g$ for some $\lambda \in \mathbb{R}$. At present no general existence results on Einstein metrics are known, except of some important classes of Einstein metrics, such as K\"{a}hler-Einstein metrics (see \cite{Ti}), Sasakian-Einstein metrics (see \cite{BGa}), homogeneous Einstein metrics (see \cite{BWZ} for the compact case and \cite{Heb} for the noncompact case). We are concerned with homogeneous Einstein metrics on reductive homogeneous spaces whose isotropy representation decomposes into a direct sum of irreducible non-equivalent summands. Then the Einstein equation reduces to a non linear algebraic system of equations. For more details on Einstein manifolds we refer to \cite{Bes}.

Nowadays, homogeneous Einstein metrics on flag manifolds have been better understood. They have been completely classified for any flag manifold $M = G/H$ of a compact simple Lie group $G$ with two (see \cite{ACh2, Sak1}), three (see \cite{Arv, Kim}) and four isotropy summands (see \cite{ACh3, ACS1, ACS2}). The case of five isotropy summands was partially solved in the work \cite{ACS3}. The case of six isotropy summands have not been completely classified except only the full flag manifold $G_2/T$ (see \cite{ACS4}) and some flag manifolds of exceptional Lie groups (see \cite{WZ}).

In the present paper we study the flag manifold $SO_{2(n_1+n_2+n_3)+1}/U_{n_1} \times U_{n_2} \times SO_{2n_3+1}$, where $n_1 > 1, \ n_2 > 1, \ n_3 > 1$, whose isotropy representation decomposes into six irreducible $Ad(H)$-submodules
$$
\mathfrak{m} = T_{eH} M = \mathfrak{m}_1 \oplus \mathfrak{m}_2 \oplus \mathfrak{m}_3 \oplus \mathfrak{m}_4 \oplus \mathfrak{m}_5 \oplus \mathfrak{m}_6.
$$
More precisely, we deal with the question about the number of complex invariant Einstein metrics on this manifold. We prove that this number is equal to exactly 132 for general parameters $(n_1, n_2, n_3)$ and provide algebraic equations on these parameters when it is not true. The proof is based on the technique of Newton polytope that firstly was applied in this setting by M. Graev (see \cite{GrUMN, GrTMMO,GrC}). Also we give new examples of Ricci-flat Lorentzian manifolds that appear as Inonu-Wigner contractions of the manifold $SO_{2(n_1+n_2+n_3)+1}/U_{n_1} \times U_{n_2} \times SO_{2n_3+1}$ by certain two-dimensional faces of the polytope.

The paper is organized as follows. In Section 2 we discuss complex $G$-invariant metrics on flag manifolds $G/H$ and give necessary formulas for metric invariants. In Section 3 we apply Bernstein-Koushnirenko theorem to Einstein equations and describe Newton polytope of the scalar curvature. Also in Section 3 we give definition of this polytope through $T$-root system of flag manifold. Section 4 is devoted to proof of the main result. At last, in Section 5 we define Inonu-Wigner contractions of homogeneous manifolds and apply this to obtain new examples of Ricci-flat Lorentzian manifolds.

\section{Invariant Einstein metrics}

Let $G / H$ be a flag manifold of a compact simple Lie group $G$. Then $H$-module of isotropy $\mathfrak{g}/\mathfrak{h}$ splits into direct sum of pairwise non-equivalent irreducible submodules. Denote by $\mathfrak{m} \subset \mathfrak{g}$ orthogonal complement of the algebra $\mathfrak{h}$ with respect to the Killing form $B$. The group $\text{Ad}_{\mathfrak{g}}(H)$ naturally acts on $\mathfrak{m}$, so that it can be identified with the module $\mathfrak{g}/\mathfrak{h}$. Assume that the module $\mathfrak{m}=\bigoplus\limits_{i=1}^{d} \mathfrak{m}_{i}$ decomposes into non-equivalent irreducible submodules $\mathfrak{m}_i, \ i=1,...,d$. Then an invariant Riemannian metric $g$ on $G/H$, its Ricci tensor $\text{Ric}$ and scalar curvature $s$ can be written as follows
\begin{equation}
g=\sum\limits_{i=1}^{d} t_i (-B|_{\mathfrak{m}_i}), \ \ \ \text{Ric}=\sum\limits_{i=1}^{d} r_i t_i (-B|_{\mathfrak{m}_i}), \ \ \ s=\sum\limits_{i=1}^{d} N_i \lambda_i, \ \ \ N_i=\text{dim}~\mathfrak{m}_i.
\end{equation}

\noindent Here $t=(t_1,...,t_d)$ belongs to the cone $\mathbb{R}_{+}^{d}$, $r_i=r_i(t)$ and $s=s(t)$ are homogeneous functions on $t$ of the degree of homogeneity $-1$. In fact, the functions $r_{i}(t)$ and $s(t)$ are Laurent polynomials on $t$ due to following formulas
\begin{equation}
r_{i}=-\frac{t_i}{N_i}\frac{\partial s}{\partial t_i}, \ \ i=1,...,d,
\end{equation}
\begin{equation}
2s(t)=\sum\limits_{i=1}^{d}\frac{N_i}{t_i}-\frac{1}{3!}\sum\limits_{i=1}^{d}\sum\limits_{j=1}^{d}\sum\limits_{k=1}^{d} {i \brack j ~ k} \frac{t_i^2+t_j^2+t_k^2}{t_i t_j t_k}.
\end{equation}

\noindent Here $b_{i,j,k} \geq 0$ are non-negative coefficients that are invariant under all permutations of indices $i,j,k$, such that
\begin{equation}
T(i,j,k): \ \ \ {i \brack j ~ k} > 0 \Longleftrightarrow B(\mathfrak{m}_i,[\mathfrak{m}_j,\mathfrak{m}_k]) \neq 0.
\end{equation}

\noindent Below is an expression of these coefficients through $(-B)$-orthonormal bases of the modules $\mathfrak{m}_i,\mathfrak{m}_j,\mathfrak{m}_k$
\begin{equation}
{i \brack j ~ k}=\sum\limits_{a=1}^{N_i}\sum\limits_{b=1}^{N_j}\sum\limits_{c=1}^{N_k} (B(X_{ai},[X_{bj},X_{ck}]))^2.
\end{equation}

\noindent The coefficients ${i \brack j ~ k}=(B(\mathfrak{m}_i),[\mathfrak{m}_j,\mathfrak{m}_k])^2$ are sometimes called the structure constants of the manifold $G/H$ and the set $T(G/H)$ of triples $(i,j,k)$ for which ${i \brack j ~ k} \neq 0$ is called the invariant of de Siebenthal.

Therefore, Einstein equations for the invariant metric $g$ can be written in the form of the following system of $d-1$ homogeneous equations for $d$ variables
\begin{equation}\label{Einstein}
r_{i}(t)-r_{i+1}(t)=0, \ \ i=1,...,d-1.
\end{equation}

\noindent Moreover, we can expand these equations to invariant complex metrics $g$ on $G/H$ which can be defined as well as Riemannian metrics just by replacing $t \in \mathbb{R}_{+}^{d}$ with $t \in (\mathbb{C}^{*})^{d}$. Thus, we can consider isolated complex solutions of the system (\ref{Einstein}) up to multiplication by complex number.

\section{Newton polytopes}

Further, we will need the Bernstein-Koushnirenko theorem, so let us remind this theorem. Suppose we have the system of algebraic equations 

In order to estimate the number of isolated complex solutions of (\ref{Einstein}) we use the compact convex polytope $\Delta=\Delta_{G/H} \subset \mathbb{R}^{d}$ that is the Newton polytope of the Laurent polynomial $s(t)$, i. e. the convex hull of finite set $\text{supp}(s)\subset \mathbb{Z}^d$ of vector exponents of all monomials of $s(t)$. Also introduce $(d-1)$-dimensional simplex $S$ with vertices $-e_i, i=1,...,d$, and define a normalized volume of polytope $\Delta$ as follows
\begin{equation}
\nu(\Delta)=\frac{\text{vol}(\Delta)}{\text{vol}(S)}=\frac{(d-1)!}{\sqrt{d}} \text{vol}(\Delta) \in \mathbb{Z}_{\geq 0}.
\end{equation}

\noindent Then straightforward applying the Bernstein-Koushnirenko theorem to Einstein equations gives us the following theorem.

\begin{Teo}[{\cite[Thm 1.7]{GrTMMO}}]
\textit{ The number $\mathcal{E}(G/H)$ of isolated complex solutions of Einstein equations (\ref{Einstein}) satisfies the inequality $\mathcal{E}(G/H) \leq \nu(\Delta)$. Moreover, the equality $\mathcal{E}(G/H)=\nu(\Delta)$ holds if and only if the hypersurface $s_{\Gamma}(t)=0$ has no singular points in $(\mathbb{C}^{*})^{d}$ for any proper face $\Gamma$ of the polytope $\Delta$. Also in case of $\mathcal{E}(G/H)=\nu(\Delta)$ all solutions of Einstein equations (\ref{Einstein}) are isolated.}
\end{Teo}
\noindent Here $s_{\Gamma}(t)$ for a face $\Gamma \subset \Delta$ denotes the sum of all monomials of $s(t)$ whose vector exponents belong to $\Gamma$ and corresponding condition on the hypersurface $s_{\Gamma}(t)=0$ we will call discriminant condition for a face $\Gamma \subset \Delta$. It is convenient to write these discriminant conditions in non-homogeneous form. More precisely, let $v \in \Delta_{\alpha}$ be a vertex of the face $\Gamma$, then after multiplication the function $s_{\Gamma}$ by monomial $t^{-v}$ we obtain the Laurent polynomial $t^{-v}s_{\Gamma}$ of the degree of homogeneity 0 essentially depending on $\text{dim}~\Gamma \leq d-2$ variables. Discriminant conditions mean that the complex hypersurface $t^{-v}s_{\Gamma}=0$ has no singular points $t$ such that $t_i \neq 0$, i. e.
\begin{equation}\label{Bernstein}
t^{-v}s_{\Gamma}=0, \ \ \ \text{d}(t^{-v}s_{\Gamma})=0, \ \ \ t_i \neq 0.
\end{equation}

It is important that some faces can be deleted from consideration due to their geometric properties because corresponding discriminant conditions necessarily hold. For example, it is obvious that for any vertex $v$ we have that $t^{-v}s_{v}(t)$ is non-zero constant, so discriminant condition holds. Further, generalization of vertex is $n$-dimensional face $\Gamma=\text{Conv}(v \cup \Gamma')$ such that all vertices except one belong to $(n-1)$-dimensional subface $\Gamma' \subset \Gamma$. We call such faces pyramidal and we have the following lemma for them.

\begin{Lem}[{\cite[Thm 5.1]{GrTMMO}}]\label{pyramid}
\textit{For pyramidal faces corresponding discriminant conditions always hold.}
\end{Lem}

Other type of face that we call octahedral for which discriminant conditions almost always hold is defined by the following way. It is a convex hull of the points $c+z_i, \ c-z_i$, where $z_i, \ i=1,...,n$ are arbitrary linear independent vectors in $\mathbb{R}^{d}$, and $c \in \mathbb{R}^{d}$ is some vector called center of octahedral face. For these faces we have the lemma.

\begin{Lem}[{\cite[Lemma 1.12]{GrTMMO}}]\label{octahedron}
\textit{Let $\Gamma$ be octahedral face with center $-e_{i_0}$ such that all other points $-e_{i}, i \neq i_0$ don't belong to $\Gamma$. Then corresponding discriminant condition holds if and only if $[\mathfrak{m}_{i_0},\mathfrak{m}_{i_0}] \neq 0$.}
\end{Lem}

Discriminant conditions for all other faces of Newton polytope $\Delta$ need to be checked separately every time. Therefore, we are mostly interested in faces that are not pyramidal and not octahedral. We will call such faces marked. Simplest example of marked faces is a parallelogram that can be defined as a convex hull of vectors $v, v+v_1, v+v_2, v+v_1+v_2$ where vectors $v_1, v_2$ are linear independent. For faces that are parallelograms the following lemma is true.
\begin{Lem}\label{prlm}
\textit{If a face $\Gamma \subset \Delta$ is a parallelogram containing only four integer points then the discriminant condition is equivalent to vanishing of determinant of matrix which elements are coefficients of $s_{\Gamma}(t)$.}
\end{Lem}

It is worth to note that Newton polytope $\Delta$ can be also described by so-called $T$-root system of a flag manifold $G/H$. For this consider characters of the center $Z(H)$ formed the lattice in $\mathbb{R}^{n}$. Denote by $\omega_i$ and $-\omega_i$ characters that are involved in the representation of torus $Z(H)$ on $\mathbb{C} \otimes \mathfrak{m}_i$. The system of all characters that are involved in the representation $Z(H)$ on $\mathbb{C} \otimes \mathfrak{m}=\bigoplus\limits_{i=1}^{d} \mathbb{C} \otimes \mathfrak{m}_i$, $\Omega=\{\omega_i,-\omega_i:i=1,...,d\}$ consists of exactly $2d$ elements, i. e. the module $\mathfrak{m}_i$ is defined by the character up to sign before $\omega_i$.

It is easy to see that $\Omega$ can be considered as $n$-dimensional system of vectors in $\mathbb{R}^{n}$, where $n=\text{dim}~Z(H)=b_2(G/H)$. This system of vectors $\Omega$ is called $T$-root system of a flag manifold $G/H$.

\begin{Teo}[{\cite[Thm 2.1]{GrTMMO}}]
\textit{Let $G/H$ be a flag manifold with $T$-root system $\Omega$. Then support $\text{supp}(s)$ of the scalar curvature $s(t)$ consists of all vertices $-e_i, i=1,...,d$, and vertices of the form $e_i-e_j-e_k$ for all $(i,j,k)$ such that the equality $\omega_i \pm \omega_j \pm \omega_k = 0$ on positive roots of $\Omega$ holds for some distribution of signs.}
\end{Teo}

\noindent As consequence we have that the Newton polytope $\Delta$ associated with a flag manifold $G/H$ depends only on its $T$-root system $\Omega$, i. e. $\Delta=\Delta(\Omega)$. Moreover, in order to classify polytopes $\Delta(\Omega)$ it is enough to consider systems of $T$-roots up to maps $\Omega \rightarrow \Omega'$ continuing to isomorphisms of their linear hulls $\mathbb{R}^{n} \rightarrow \mathbb{R}^{n'}$. Also note that the dimension of the Newton polytope is equal to $d-1$.

\section{Flag manifold $M_{n_1,n_2,n_3}=SO_{2(n_1+n_2+n_3)+1}/U_{n_1} \times U_{n_2} \times SO_{2n_3+1}$}

Note that flag manifold $SO_{2(n_1+n_2+n_3)+1}/U_{n_1} \times U_{n_2} \times SO_{2n_3+1}$, where $n_1>1, \ n_2>1, \ n_3 > 1,$ has the $T$-root system $BC_2$ that can be represented by the following painted Dynkin diagram

\begin{center}\includegraphics[scale=0.44]{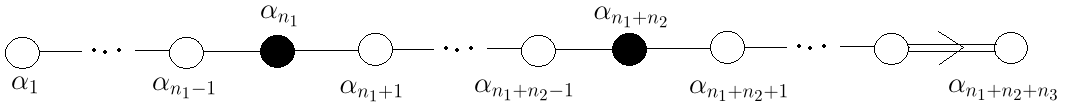}\end{center}

\noindent Choose simple roots $\overline{\alpha_{n_1}}, \overline{\alpha_{n_1+n_2}}$ of $BC_2$ such that the following roots are positive
$$
\omega_1=2 \overline{\alpha_{n_1}} + 2 \overline{\alpha_{n_1 + n_2}}, \ \ \omega_2=2 \overline{\alpha_{n_1 + n_2}}, \ \ \omega_3=\overline{\alpha_{n_1}},
$$
$$
\omega_4=\overline{\alpha_{n_1}}+2 \overline{\alpha_{n_1 + n_2}}, \ \ \omega_5 = \overline{\alpha_{n_1}} + \overline{\alpha_{n_1 + n_2}}, \ \ \omega_6 = \overline{\alpha_{n_1 + n_2}}.
$$


\noindent It is obvious that these positive roots $\omega_i$ satisfy the following relations
$$\omega_1-\omega_3-\omega_4=0,~~~~~~~~~~~~  \omega_3-\omega_5+\omega_6=0,~~~~~~~~~~~~ \omega_1-2\omega_5=0, $$
$$\omega_2+\omega_3-\omega_4=0,~~~~~~~~~~~~ \omega_4-\omega_5-\omega_6=0,~~~~~~~~~~~~ \omega_2-2\omega_6=0.$$
\noindent Hence, the invariant of de Siebenthal of the flag manifold $M_{n_1, n_2, n_3}$ is the following
$$T(M_{n_1, n_2, n_3})=\{\text{all permutations of triples}~ (1,3,4), (2,3,4), (3,5,6), (4,5,6), (1,5,5), (2,6,6)\}.$$
Therefore, it is necessary to compute only structure constants ${1 \brack 3 ~ 4}, {2 \brack 3 ~ 4}, {3 \brack 5 ~ 6}, {4 \brack 5 ~ 6}, {1 \brack 5 ~ 5}, {2 \brack 6 ~ 6}$.

Fundamental weights of the root system $B_{n_1+n_2+n_3}$ are the following
$$
\Lambda_{1}=\alpha_{1}+...+\alpha_{n_1+n_2+n_3},
$$
$$
\Lambda_{2}=\alpha_{1}+ 2 \big( \alpha_{2} +...+\alpha_{n_1+n_2+n_3} \big),
$$
$$
\Lambda_{3}=\alpha_{1}+ 2 \alpha_{2} + 3 \big( \alpha_{3} + ...+\alpha_{n_1+n_2+n_3} \big),
$$
$$
\Lambda_{4}=\alpha_{1}+ 2 \alpha_{2} + 3 \alpha_{3} + 4 \big( \alpha_{4} + ... +\alpha_{n_1+n_2+n_3} \big),
$$
{\hspace{40mm}.....................................................................}
$$
\Lambda_{n_1+n_2+n_3}=\frac{1}{2} \big( \alpha_{1}+2\alpha_{2}+3\alpha_{3}+...+(n_1+n_2+n_3) \alpha_{n_1+n_2+n_3} \big).
$$

Fix orthogonal decomposition $\mathfrak{g}=\mathfrak{h} \oplus \mathfrak{m}$ with respect to the Killing form. Then, according to the first section, $H$-module $\mathfrak{m}$ splits into direct sum of 6 pairwise inequivalent irreducible submodules $\mathfrak{m}_i, i=1,...,6$ corresponding to positive roots $\omega_i$.

\begin{Teo}
\textit{The submodules $m_{i}$ have the following dimensions}
$$N_1=\emph{dim}~\mathfrak{m}_{1}=n_1(n_1-1),\ \ \ N_2=\emph{dim}~\mathfrak{m}_{2}=n_2(n_2-1),$$
$$N_3=\emph{dim}~\mathfrak{m}_{3}=2 n_1 n_2,\ \ \ N_4=\emph{dim}~\mathfrak{m}_{4}=2 n_1 n_2,$$
$$N_5=\emph{dim}~\mathfrak{m}_{5}=2n_1(2n_3+1),\ \ \ N_6=\emph{dim}~\mathfrak{m}_{6}=2n_2(2n_3+1).$$
\end{Teo}
\noindent \textit{Proof:} Use the following Weyl formula
$$
\text{dim}_{\mathbb{C}}~\mathfrak{m}_{i}=\prod_{\phi \in R_{K}^{+}} \Bigg( 1+\frac{(\lambda_{i},\phi)}{(\delta_{K},\phi)} \Bigg)
$$
$$
\lambda_{1}=\alpha_1+...+\alpha_{n_1+n_2-1}=\Lambda_{1}+\Lambda_{n_1+n_2-1}-\Lambda_{n_1+n_2},
$$
$$
\lambda_{2}=\alpha_{n_1+1}+...+\alpha_{n_1+n_2}+2(\alpha_{n_1+n_2+1}+...+\alpha_{n_1+n_2+n_3})=\Lambda_{n_1+1}-\Lambda_{n_1}+\Lambda_{n_1+n_2+1}-\Lambda_{n_1+n_2},
$$
$$
\lambda_{3}=\alpha_{1}+...+\alpha_{n_1+n_2}+2(\alpha_{n_1+n_2+1}+...+\alpha_{n_1+n_2+n_3})=\Lambda_{1}+\Lambda_{n_1+n_2+1}-\Lambda_{n_1+n_2},
$$
$$
\lambda_{4}=\alpha_{n_1+1}+2(\alpha_{n_1+2}+...+\alpha_{n_1+n_2+n_3})=\Lambda_{n_1+2}-\Lambda_{n_1},
$$
$$
\lambda_{5}=\alpha_1+...+\alpha_{n_1}+2(\alpha_{n_1+1}+ ... + \alpha_{n_1+n_2+n_3})=\Lambda_{1}+\Lambda_{n_1+1}-\Lambda_{n_1},
$$
$$
\lambda_{6}=\alpha_1+ 2 (\alpha_2 + ... + \alpha_{n_1+n_2+n_3})=\Lambda_{2},
$$
$$
\delta_{K}=\Lambda_{1}+...+\Lambda_{n_1-1}+\Lambda_{n_1+1}+...+\Lambda_{n_1+n_2-1}+\Lambda_{n_1+n_2+1}+...\Lambda_{n_1+n_2+n_3}
$$
$$
\text{dim}_{\mathbb{C}}~\mathfrak{m}_{1}=\Bigg( 1 + \frac{1}{1} \Bigg)...\Bigg( 1+\frac{1}{n_1 - 2} \Bigg) \cdot \Bigg( 1 + \frac{1}{2} \Bigg)...\Bigg( 1+\frac{1}{n_1 - 1} \Bigg)=\frac{n_1 (n_1 - 1)}{2}.
$$
$$
\text{dim}_{\mathbb{C}}~\mathfrak{m}_{2}=\Bigg( 1 + \frac{1}{1} \Bigg)...\Bigg( 1+\frac{1}{n_2 - 2} \Bigg) \cdot \Bigg( 1 + \frac{1}{2} \Bigg)...\Bigg( 1+\frac{1}{n_2 - 1} \Bigg)=\frac{n_2 (n_2 - 1)}{2},
$$
$$
\text{dim}_{\mathbb{C}}~\mathfrak{m}_{3}=\Bigg( 1+\frac{1}{1} \Bigg)...\Bigg( 1+\frac{1}{n_1-1} \Bigg) \cdot \Bigg( 1+\frac{1}{1} \Bigg)...\Bigg( 1+\frac{1}{n_2-1} \Bigg)=n_1 n_2,
$$
$$
\text{dim}_{\mathbb{C}}~\mathfrak{m}_{4}=\Bigg( 1+\frac{1}{1} \Bigg)...\Bigg( 1+\frac{1}{n_1-1} \Bigg) \cdot \Bigg( 1+\frac{1}{1} \Bigg)...\Bigg( 1+\frac{1}{n_2-1} \Bigg)=n_1 n_2,
$$
$$
\text{dim}_{\mathbb{C}}~\mathfrak{m}_{5}=\Bigg( 1+\frac{1}{1} \Bigg)...\Bigg( 1+\frac{1}{n_1-1} \Bigg) \cdot \Bigg( 1+\frac{1}{1} \Bigg) ... \Bigg( 1+\frac{1}{n_3-1} \Bigg) \Bigg( 1+\frac{2}{2(n_3-1)+1} \Bigg) \cdot
$$
$$
\cdot \Bigg( 1 + \frac{1}{1+2(n_3-2)+1} \Bigg) ... \Bigg( 1 + \frac{1}{(n_3-1)+ 2 \cdot 0 + 1 \cdot 1} \Bigg) = n_1 (2 n_3 + 1),
$$
$$
\text{dim}_{\mathbb{C}}~\mathfrak{m}_{6}=\Bigg( 1+\frac{1}{1} \Bigg)...\Bigg( 1+\frac{1}{n_2-1} \Bigg) \cdot \Bigg( 1+\frac{1}{1} \Bigg) ... \Bigg( 1+\frac{1}{n_3-1} \Bigg) \Bigg( 1+\frac{2}{2(n_3-1)+1} \Bigg) \cdot
$$
$$
\cdot \Bigg( 1 + \frac{1}{1+2(n_3-2)+1} \Bigg) ... \Bigg( 1 + \frac{1}{(n_3-1)+ 2 \cdot 0 + 1 \cdot 1} \Bigg) = n_2 (2 n_3 + 1),
$$
\hfill$\Box$

\noindent\textbf{Remark}: It can be shown also just by dealing with lengths of white strings.

\noindent Hence, the polynomial of scalar curvature $s(t)$ and components of Ricci tensor are defined by the following formulas
$$
2s(t)=N_1 t_1^{-1}+N_2 t_2^{-1}+N_3 t_3^{-1}+N_4 t_4^{-1}+N_5 t_5^{-1}+N_6 t_6^{-1} -
$$
$$
- {1 \brack 3 ~ 4}\Big(t_1 t_3^{-1} t_4^{-1} + t_1^{-1} t_3 t_4^{-1} + t_1^{-1} t_3^{-1} t_4\Big) - {2 \brack 3 ~ 4} \Big( t_2 t_3^{-1} t_4^{-1} + t_2^{-1} t_3 t_4^{-1} + t_2^{-1} t_3^{-1} t_4 \Big) -$$$$- {3 \brack 5 ~ 6} \Big( t_3 t_5^{-1} t_6^{-1} + t_3^{-1} t_5 t_6^{-1} + t_3^{-1} t_5^{-1} t_6 \Big) - {4 \brack 5 ~ 6} \Big( t_4 t_5^{-1} t_6^{-1} + t_4^{-1} t_5 t_6^{-1} + t_4^{-1} t_5^{-1} t_6 \Big) -$$$$- \frac{1}{2} {1 \brack 5 ~ 5} \Big ( 2 t_1^{-1} + t_1 t_5^{-2} \Big) - \frac{1}{2} {2 \brack 6 ~ 6} \Big( 2 t_2^{-1} + t_2 t_6^{-2} \Big),
$$
\begin{numcases}{}
r_1=\frac{t_1^{-1}}{2} + \frac{1}{2 N_1} {1 \brack 3 ~ 4} \Big( t_1 t_3^{-1} t_4^{-1} - t_1^{-1} t_3 t_4^{-1} - t_1^{-1} t_3^{-1} t_4 \Big) + \frac{1}{4 N_1} {1 \brack 5 ~ 5} \Big( -2 t_1^{-1} + t_1 t_5^{-2} \Big), \nonumber \\
r_2=\frac{t_2^{-1}}{2} + \frac{1}{2 N_2} {2 \brack 3 ~ 4} \Big( t_2 t_3^{-1} t_4^{-1} - t_2^{-1} t_3 t_4^{-1} - t_2^{-1} t_3^{-1} t_4 \Big) + \frac{1}{4 N_2} {2 \brack 6 ~ 6} \Big( -2 t_2^{-1} + t_2 t_6^{-2} \Big), \nonumber \\
r_3=\frac{t_3^{-1}}{2} + \frac{1}{2 N_3} {1 \brack 3 ~ 4} \Big( -t_1 t_3^{-1} t_4^{-1} + t_1^{-1} t_3 t_4^{-1} - t_1^{-1} t_3^{-1} t_4 \Big) +  \nonumber \\
+ \frac{1}{2 N_3} {2 \brack 3 ~ 4} \Big( - t_2 t_3^{-1} t_4^{-1} + t_2^{-1} t_3 t_4^{-1} - t_2^{-1} t_3^{-1} t_4 \Big) + \frac{1}{2 N_3} {3 \brack 5 ~ 6} \Big( t_3 t_5^{-1} t_6^{-1} - t_3^{-1} t_5 t_6^{-1} - t_3^{-1} t_5^{-1} t_6 \Big), \nonumber \\
r_4=\frac{t_4^{-1}}{2} + \frac{1}{2 N_4} {1 \brack 3 ~ 4} \Big( -t_1 t_3^{-1} t_4^{-1} - t_1^{-1} t_3 t_4^{-1} + t_1^{-1} t_3^{-1} t_4 \Big) +  \nonumber \\
+ \frac{1}{2 N_4} {2 \brack 3 ~ 4} \Big( - t_2 t_3^{-1} t_4^{-1} - t_2^{-1} t_3 t_4^{-1} + t_2^{-1} t_3^{-1} t_4 \Big) + \frac{1}{2 N_4} {4 \brack 5 ~ 6} \Big( t_4 t_5^{-1} t_6^{-1} - t_4^{-1} t_5 t_6^{-1} - t_4^{-1} t_5^{-1} t_6 \Big), \nonumber \\
r_5=\frac{t_5^{-1}}{2} + \frac{1}{2 N_5} {3 \brack 5 ~ 6} \Big( - t_3 t_5^{-1} t_6^{-1} + t_3^{-1} t_5 t_6^{-1} - t_3^{-1} t_5^{-1} t_6 \Big) + \nonumber \\
+ \frac{1}{2 N_5} {4 \brack 5 ~ 6} \Big( - t_4 t_5^{-1} t_6^{-1} + t_4^{-1} t_5 t_6^{-1} - t_4^{-1} t_5^{-1} t_6 \Big) - \frac{1}{2 N_5} {1 \brack 5 ~ 5}  t_1 t_5^{-2}, \nonumber \\
r_6=\frac{t_6^{-1}}{2} + \frac{1}{2 N_6} {3 \brack 5 ~ 6} \Big( - t_3 t_5^{-1} t_6^{-1} - t_3^{-1} t_5 t_6^{-1} + t_3^{-1} t_5^{-1} t_6 \Big) + \nonumber \\
+ \frac{1}{2 N_6} {4 \brack 5 ~ 6} \Big( - t_4 t_5^{-1} t_6^{-1} - t_4^{-1} t_5 t_6^{-1} + t_4^{-1} t_5^{-1} t_6 \Big) - \frac{1}{2 N_6} {2 \brack 6 ~ 6}  t_2 t_6^{-2}.
\end{numcases}

\begin{Teo}\label{constants}
\textit{The structure constants of the flag manifold $M_{n_1, n_2, n_3}=SO_{2(n_1+n_2+n_3)+1}/U_{n_1} \times U_{n_2} \times SO_{2n_3+1}$, where $n_1>1, \ n_2>1, \ n_3>1$, are the following}
$${1 \brack 3 ~ 4}=\frac{n_1 n_2(n_1-1)}{2(n_1+n_2+n_3)-1},~~~~ {2 \brack 3 ~ 4}=\frac{n_1 n_2 (n_2-1)}{2(n_1+n_2+n_3)-1},$$
$${3 \brack 5 ~ 6}={4 \brack 5 ~ 6}=\frac{n_1 n_2 (2 n_3 + 1)}{2(n_1+n_2+n_3)-1},$$
$${1 \brack 5 ~ 5}=\frac{n_1(n_1-1)(2n_3+1)}{2(n_1+n_2+n_3)-1},~~~~ {2 \brack 6 ~ 6}=\frac{n_2(n_2-1)(2n_3+1)}{2(n_1+n_2+n_3)-1}.$$
\end{Teo}

\noindent \textit{Proof:} In order to compute the structure constants ${i \brack j ~ k}$ we consider the following Riemannian submersion
$$
H/L \hookrightarrow G/H \twoheadrightarrow G/L,
$$
\noindent where
$$
H \hookrightarrow L \hookrightarrow G,
$$
$$
H=U_{n_1} \times U_{n_2} \times SO_{2n_3+1}, \ L=U_{n_1+n_2} \times SO_{2n_3+1}, \ G=SO_{2(n_1+n_2+n_3)+1}.
$$
\noindent Corresponding decomposition of the Lie algebra $\mathfrak{g}=\mathfrak{s}\mathfrak{o}_{2(n_1+n_2+n_3)+1}$ has the form
$$
\mathfrak{g}=\mathfrak{h} \oplus \mathfrak{q} \oplus \mathfrak{p}, \ \ \text{where} \ \ \text{Lie}(L)=\mathfrak{h} \oplus \mathfrak{q},
$$
$$
\mathfrak{q}=\mathfrak{m}_3 \ \ \text{and} \ \ \mathfrak{p}=\mathfrak{m}_1 \oplus \mathfrak{m}_2 \oplus \mathfrak{m}_4 \oplus \mathfrak{m}_5 \oplus \mathfrak{m}_6.
$$
\noindent Next, the $L$-module of isotropy $\mathfrak{g}/(\mathfrak{\mathfrak{h} \oplus \mathfrak{q}}) \simeq \mathfrak{p}$ has two non-equivalent irreducible submodules $\mathfrak{p}_1=\mathfrak{m}_5 \oplus \mathfrak{m}_6$ and $\mathfrak{p}_2=\mathfrak{m}_1 \oplus \mathfrak{m}_2 \oplus \mathfrak{m}_4$, respectively. The scalar curvature of $L$-invariant metric $\check{g}=x_1 (-B)|_{\mathfrak{p}_1} \oplus x_2 (-B)|_{\mathfrak{p}_2}$ on $G/L$ is equal to
$$
s(\check{g})=\frac{d_1}{2 x_1} + \frac{d_2}{2 x_2} - \frac{1}{4}\frac{d_1 d_2}{d_1 + 4 d_2} \Bigg( \frac{2}{x_2} + \frac{x_2}{x^{2}_{1}} \Bigg),
$$
\noindent where
$$
d_{1}=\text{dim~} \mathfrak{p}_1=N_{5}+N_{6}, \ \ \ d_{2}=\text{dim~} \mathfrak{p}_2=N_1+N_2+N_4.
$$
\noindent so the components $\check{r}_{1}$ and $\check{r}_{2}$ of the corresponding Ricci tensor have the following form
$$
\check{r}_{1}=\frac{1}{2 x_1} - \frac{1}{2} \frac{d_2}{d_1 + 4 d_2} \frac{x_2}{x^{2}_{1}}, \ \ \ \check{r}_{2}=\frac{1}{2 x_2} + \frac{1}{4} \frac{d_{1}}{d_{1} + 4 d_{2}} \Bigg( - \frac{2}{x_2} + \frac{x_2}{x^{2}_{1}} \Bigg).
$$
On the other hand, the horizontal part of $\text{Ric}(g)$ coincides with $\text{Ric}(\check{g})$ (see [?]). It means that if we compute the components $r_i$ of Ricci tensor $Ric(g)$ for the $H$-invariant metric $g=(x_{2}, x_{2}, y, x_{2}, x_{1}, x_{1})$, then the subexpressions of $r_5, r_6$ depending only on $x_1, x_2$ coincide with $\check{r}_{1}$, and the ones of $r_1, r_2, r_4$ with $\check{r}_{2}$, respectively. More precisely, the components of $\text{Ric}(g)$ for such metric are the following
\begin{numcases}{}
r_1=\frac{x_2^{-1}}{2} - \frac{1}{2 N_1} {1 \brack 3 ~ 4} x_2^{-2} y + \frac{1}{4 N_1} {1 \brack 5 ~ 5} \Big( -2 x_2^{-1} + x_2 x_1^{-2} \Big), \nonumber \\
r_2=\frac{x_2^{-1}}{2} - \frac{1}{2 N_2} {2 \brack 3 ~ 4} x_2^{-2} y + \frac{1}{4 N_2} {2 \brack 6 ~ 6} \Big( -2 x_2^{-1} + x_2 x_1^{-2} \Big), \nonumber \\
r_3=\frac{y^{-1}}{2} + \frac{1}{2 N_3} \Bigg( {1 \brack 3 ~ 4} + {2 \brack 3 ~ 4} \Bigg) \Big( - 2 y^{-1} + x_2^{-2} y \Big) +\frac{1}{2 N_3} {3 \brack 5 ~ 6} \Big( y x_1^{-2} - 2 y^{-1} \Big), \nonumber \\
r_4=\frac{x_2^{-1}}{2} - \frac{1}{2 N_4} \Bigg( {1 \brack 3 ~ 4} + {2 \brack 3 ~ 4} \Bigg) x_2^{-2} y + \frac{1}{2 N_4} {4 \brack 5 ~ 6} \Big( x_2 x_1^{-2} - 2 x_2^{-1} \Big), \nonumber \\
r_5=\frac{x_1^{-1}}{2} - \frac{1}{2 N_5} {3 \brack 5 ~ 6} y x_1^{-2} - \frac{1}{2 N_5} \Bigg( {4 \brack 5 ~ 6} + {1 \brack 5 ~ 5} \Bigg) x_2 x_1^{-2}, \nonumber \\
r_6=\frac{x_1^{-1}}{2} - \frac{1}{2 N_6} {3 \brack 5 ~ 6} y x_1^{-2} - \frac{1}{2 N_6} \Bigg( {4 \brack 5 ~ 6} + {2 \brack 6 ~ 6} \Bigg) x_2 x_1^{-2}.
\end{numcases}
\noindent So we obtain the following relations
$$
\frac{1}{4 N_1} {1 \brack 5 ~ 5} =  \frac{1}{4 N_2} {2 \brack 6 ~ 6} = \frac{1}{2 N_4} {4 \brack 5 ~ 6} = \frac{1}{4} \frac{d_{1}}{d_{1} + 4 d_{2}}=\frac{1}{4} \frac{N_5+N_6}{N_5+N_6 + 4 (N_1+N_2+N_4)},
$$
$$
\frac{1}{2 N_5} \Bigg( {4 \brack 5 ~ 6} + {1 \brack 5 ~ 5} \Bigg)=\frac{1}{2 N_6} \Bigg( {4 \brack 5 ~ 6} + {2 \brack 6 ~ 6} \Bigg)=\frac{1}{2} \frac{d_2}{d_1 + 4 d_2}=\frac{1}{2} \frac{N_1+N_2+N_4}{N_5+N_6 + 4 (N_1+N_2+N_4)}.
$$
\noindent From these relations we immediately find out several structure constants
$$
{1 \brack 5 ~ 5}=\frac{N_1(N_5+N_6)}{N_5+N_6 + 4 (N_1+N_2+N_4)}=\frac{n_1(n_1-1)(2n_3+1)}{2(n_1+n_2+n_3)-1},
$$
$$
{2 \brack 6 ~ 6}=\frac{N_2(N_5+N_6)}{N_5+N_6 + 4 (N_1+N_2+N_4)}=\frac{n_1(n_1-1)(2n_3+1)}{2(n_1+n_2+n_3)-1},
$$
$$
{4 \brack 5 ~ 6}=\frac{1}{2} \frac{N_4(N_5+N_6)}{N_5+N_6 + 4 (N_1+N_2+N_4)}=\frac{n_1 n_2 (2 n_3 + 1)}{2(n_1+n_2+n_3)-1}.
$$
\noindent Moreover, the symmetry given by the action of Weyl group implies that
$$
{3 \brack 5 ~ 6}=\frac{1}{2} \frac{N_3(N_5+N_6)}{N_5+N_6 + 4 (N_1+N_2+N_3)}=\frac{n_1 n_2 (2 n_3 + 1)}{2(n_1+n_2+n_3)-1}.
$$

Now, in the similar way, consider the Riemannian submersion
$$
H/L' \hookrightarrow G/H \twoheadrightarrow G/L',
$$
\noindent where
$$
H \hookrightarrow L' \hookrightarrow G,
$$
$$
H=U_{n_1} \times U_{n_2} \times SO_{2n_3+1}, \ L'=U_{n_1} \times SO_{2(n_2+n_3)+1}, \ G=SO_{2(n_1+n_2+n_3)+1}.
$$
\noindent Corresponding decomposition of the Lie algebra $\mathfrak{g}=\mathfrak{s}\mathfrak{o}_{2(n_1+n_2+n_3)+1}$ has the form
$$
\mathfrak{g}=\mathfrak{h} \oplus \mathfrak{q}' \oplus \mathfrak{p}', \ \ \text{where} \ \ \text{Lie}(L')=\mathfrak{h} \oplus \mathfrak{q}',
$$
$$
\mathfrak{q}'=\mathfrak{m}_2 \oplus \mathfrak{m}_6 \ \ \text{and} \ \ \mathfrak{p}'=\mathfrak{m}_1 \oplus \mathfrak{m}_3 \oplus \mathfrak{m}_4 \oplus \mathfrak{m}_5.
$$
\noindent Next, the $L$-module of isotropy $\mathfrak{g}/(\mathfrak{\mathfrak{h} \oplus \mathfrak{q}}) \simeq \mathfrak{p}$ has two non-equivalent irreducible submodules $\mathfrak{p}'_1=\mathfrak{m}_3 \oplus \mathfrak{m}_4 \oplus \mathfrak{m}_5$ and $\mathfrak{p}'_2=\mathfrak{m}_1$, respectively. The scalar curvature of $L$-invariant metric $\check{g}=x_1 (-B)|_{\mathfrak{p}_1} \oplus x_2 (-B)|_{\mathfrak{p}_2}$ on $G/L$ is equal to
$$
s(\check{g})=\frac{d'_1}{2 x_1} + \frac{d'_2}{2 x_2} - \frac{1}{4}\frac{d'_1 d'_2}{d'_1 + 4 d'_2} \Bigg( \frac{2}{x_2} + \frac{x_2}{x^{2}_{1}} \Bigg),
$$
\noindent where
$$
d'_{1}=\text{dim~} \mathfrak{p}'_1=N_{3}+N_{4}+N_{5}, \ \ \ d'_{2}=\text{dim~} \mathfrak{p}'_2=N_1.
$$
\noindent so the components $\check{r}_{1}$ and $\check{r}_{2}$ of the corresponding Ricci tensor have the following form
$$
\check{r}_{1}=\frac{1}{2 x_1} - \frac{1}{2} \frac{d'_2}{d'_1 + 4 d'_2} \frac{x_2}{x^{2}_{1}}, \ \ \ \check{r}_{2}=\frac{1}{2 x_2} + \frac{1}{4} \frac{d'_{1}}{d'_{1} + 4 d'_{2}} \Bigg( - \frac{2}{x_2} + \frac{x_2}{x^{2}_{1}} \Bigg).
$$
As previously if we compute the components $r_i$ of Ricci tensor $Ric(g)$ for the $H$-invariant metric $g=(x_{2}, y, x_{1}, x_{1}, x_{1}, y)$, then the subexpressions of $r_3, r_4, r_5$ depending only on $x_1, x_2$ coincide with $\check{r}_{1}$, and the ones of $r_1$ with $\check{r}_{2}$, respectively. The components of $\text{Ric}(g)$ for such metric are the following
\begin{numcases}{}
r_1=\frac{x_2^{-1}}{2} + \frac{1}{4 N_1} \Bigg( 2 {1 \brack 3 ~ 4} + {1 \brack 5 ~ 5} \Bigg) \Big( - 2 x_2^{-1} + x_2 x_1^{-2}  \Big), \nonumber \\
r_2=\frac{y^{-1}}{2} + \frac{1}{2 N_2} {2 \brack 3 ~ 4} \Big( y x_1^{-2} - 2 y^{-1} \Big) - \frac{1}{4 N_2} {2 \brack 6 ~ 6} y^{-1} , \nonumber \\
r_3=\frac{x_1^{-1}}{2} - \frac{1}{2 N_3} {1 \brack 3 ~ 4} x_2 x_1^{-2}  - \frac{1}{2 N_3} \Bigg( {2 \brack 3 ~ 4} + {3 \brack 5 ~ 6} \Bigg) y x_1^{-2}, \nonumber \\
r_4=\frac{x_1^{-1}}{2} - \frac{1}{2 N_4} {1 \brack 3 ~ 4} x_2 x_1^{-2} - \frac{1}{2 N_4} \Bigg( {2 \brack 3 ~ 4} + {4 \brack 5 ~ 6} \Bigg) y x_1^{-2}, \nonumber \\
r_5=\frac{x_1^{-1}}{2} - \frac{1}{2 N_5} \Bigg( {3 \brack 5 ~ 6} +{4 \brack 5 ~ 6} \Bigg) y x_1^{-2} - \frac{1}{2 N_5} {1 \brack 5 ~ 5}  x_2 x_1^{-2}, \nonumber \\
r_6=\frac{y^{-1}}{2} + \frac{1}{2 N_6} \Bigg( {3 \brack 5 ~ 6} + {4 \brack 5 ~ 6} \Bigg) \Big( - 2 y^{-1} + y x_1^{-2} \Big) - \frac{1}{2 N_6} {2 \brack 6 ~ 6}  y^{-1}.
\end{numcases}
\noindent So we obtain the following relations
$$
\frac{1}{N_1} \Bigg( 2 {1 \brack 3 ~ 4} + {1 \brack 5 ~ 5} \Bigg)=\frac{d'_{1}}{d'_{1} + 4 d'_{2}}=\frac{N_{3}+N_{4}+N_{5}}{N_{3}+N_{4}+N_{5} + 4 N_1},
$$
$$
\frac{1}{N_3} {1 \brack 3 ~ 4}=\frac{1}{N_5} {1 \brack 5 ~ 5}= \frac{d'_2}{d'_1 + 4 d'_2}= \frac{N_1}{N_{3}+N_{4}+N_{5} + 4 N_1}.
$$
\noindent From these relations we immediately find out remaining structure constants
$$
{1 \brack 3 ~ 4}= \frac{N_1 N_3}{N_{3}+N_{4}+N_{5} + 4 N_1}=\frac{n_1 n_2(n_1-1)}{2(n_1+n_2+n_3)-1},
$$
$$
{1 \brack 5 ~ 5}= \frac{N_1 N_5}{N_{3}+N_{4}+N_{5} + 4 N_1}.
$$
\noindent Moreover, the symmetry given by the action of Weyl group implies that
$$
{2 \brack 3 ~ 4}= \frac{N_2 N_3}{N_{3}+N_{4}+N_{6} + 4 N_2}=\frac{n_1 n_2 (n_2-1)}{2(n_1+n_2+n_3)-1},
$$
$$
{2 \brack 6 ~ 6}= \frac{N_2 N_6}{N_{3}+N_{4}+N_{6} + 4 N_2}.
$$

\hfill$\Box$

\begin{Teo}\label{number}
\textit{The number of complex invariant Einstein metrics up to homotheties on the flag manifold $M_{n_1, n_2, n_3}=SO_{2(n_1+n_2+n_3)+1}/U_{n_1} \times U_{n_2} \times SO_{2n_3+1}$, where $n_1>1, \ n_2>1, \ n_3>1$, is equal to 132, except cases when parameters $n_1, n_2, n_3$ satisfy at least one of the following algebraic equations:}
\begin{equation}\label{eqth1}
n_1 - n_2 = 0,
\end{equation}
\begin{equation}\label{eqth2}
8 n_i (2 n_3 + 1) - (n_j - 1)^2 = 0,
\end{equation}
\begin{equation}\label{eqth7}
16 n_j (n_j-1)^2 (n_i-n_j) (2 n_3 + 1) \Big(2(n_1+n_2+n_3)-1 \Big)^2 \Big(2(2n_1+2n_2+n_3)-3\Big) =
\end{equation}
$$
\Bigg(4 n_j (2 n_3 + 1) \Big( 2 (n_1 + n_2 + n_3) - 1 \Big)^2 + (n_j-1)^2 \Big( 2(2n_1 + 2n_2+ n_3) -3 \Big)(n_i - n_j) -
$$
$$
- 8 n_i^2 (2n_3+1) (4 n_i + 5 n_j + 2 n_3 - 3)\Bigg)^2,
$$
\begin{equation}\label{eqth6}
\emph{Res}(p_1,p_2)=0, \ \ \ \emph{Res}(q_{1}^{-}, q_2)=0, \ \ \ \emph{Res}(q_{1}^{+}, q_2)=0,
\end{equation}
\noindent \textit{where $(i,j)=(1,2), (2,1)$ and $p_1, p_2, q_{1}^{\mp}, q_2$ are certain univariate polynomials with coefficients depending on parameters $n_1, n_2, n_3$.}
\end{Teo}

\noindent \textit{Proof:} In order to simplify formulas, we consider the function $\tilde{s}(t)=4 \Big( 2(n_1+n_2+n_3)-1 \Big) s(t)$ which has the following expression
$$\tilde{s}(t)=a_1 t_1^{-1}+a_2 t_2^{-1}+a_3 t_3^{-1} + a_3 t_4^{-1} + a_5 t_5^{-1}+a_6 t_6^{-1} - $$$$ - b_1 \Big( t_1 t_3^{-1} t_4^{-1} + t_1^{-1} t_3 t_4^{-1} + t_1^{-1} t_3^{-1} t_4 \Big) - b_2 \Big( t_2 t_3^{-1} t_4^{-1} + t_2^{-1} t_3 t_4^{-1} + t_2^{-1} t_3^{-1} t_4 \Big) - $$$$ - b_3 \Big( t_3 t_5^{-1} t_6^{-1} + t_3^{-1} t_5 t_6^{-1} + t_3^{-1} t_5^{-1} t_6 + t_4 t_5^{-1} t_6^{-1} + t_4^{-1} t_5 t_6^{-1} + t_4^{-1} t_5^{-1} t_6 \Big) - b_4 t_1 t_5^{-2} - b_5 t_2 t_6^{-2},$$
\noindent where $a_1, ..., b_5$ are the following constants
$$a_1=4 n_1 (n_1-1) (n_1 + n_2 - 1),~~~~ a_2=4 n_2 (n_2-1) ( n_1 + n_2 - 1),$$
$$a_3=4 n_1 n_2 (2(n_1+n_2+n_3)-1),$$
$$a_4=4 n_1 (2 n_3 + 1) (2(n_1+n_2+n_3)-1),~~~~ a_5=4 n_2 (2 n_3 + 1) (2(n_1+n_2+n_3)-1),$$
$$b_1=2 n_1 n_2 (n_1-1),~~~~ b_2=2 n_1 n_2 (n_2-1),$$
$$b_3=2 n_1 n_2 (2n_3+1),$$
$$b_4=n_1 (n_1-1) (2 n_3 + 1),~~~~ b_5=n_2 (n_2-1) (2 n_3 + 1).$$

Now we have to consider Newton polytope $\Delta$ of this function or Newton polytope $\Delta(BC_2)$ of $T$-root system $BC_2$ that is the same and find out conditions when scalar curvature restricted on every face $\Gamma \subset \Delta$ has critical points in algebraic torus $(\mathbb{C}^{*})^{6}$. In fact, this polytope $\Delta$ is well described by M. Graev in the paper \cite{GrTMMO}. More precisely, it is known that $\nu(\Delta)=132$, so we have the following estimation on the number $\mathcal{E}(M_{n_1, n_2, n_3})$ of isolated Einstein metrics
\begin{equation}\label{estimate}
\mathcal{E}(M_{n_1, n_2, n_3}) \leq 132.
\end{equation}

\noindent Moreover, the marked faces of this polytope are listed in the table below where we use the following notation. Any face $\Gamma$ of the polytope $\Delta$ can be determined by the normal vector $f_{\Gamma}$ that is orthogonal to it and belongs to dual cone $\mathcal{C}(\Delta)$. If the face $\Gamma$ is facet, then vector $f_{\Gamma}$ generates the edge of dual cone. If the face $\Gamma$ is not a facet, then we choose vector $f_{\Gamma}$ that is proportional to sum of all vectors $f_{\Gamma'}$, where $\Gamma'$ are facets intersecting the face $\Gamma$. It is worth to note that the group $\text{Aut}(BC_2)$ acts on $\mathbb{R}^{6}$ by compositions of permutations of coordinates $(x_3 \leftrightarrow x_4)$ and $(x_1 \leftrightarrow x_2, \ x_5 \leftrightarrow x_6)$, so that this action continues to the set of faces of polytope $\Delta$. Since all faces from an orbit of this action are similar, only one representative of each orbit is given in the table.
\begin{center}
\begin{tabular}{|c|c|c|c|}
\hline
\multicolumn{4}{|c|}{Marked faces of the polytope $\Delta=\Delta(BC_2)$} \\
\hline
Face & Normal vector & Dimension & The number of elements of orbit \\
\hline
$\Gamma_{1}$ & $(4,2,2,4,3,1)$  & 2 & 4 \\
\hline
$\Gamma_{2}$ & $(0,0,1,1,1,1)$ & 2 & 1 \\
\hline
$\Gamma_{3}$ & $(2,0,3,3,2,1)$  & 2 & 2 \\
\hline
$\Gamma_{4}$ & $(2,0,2,2,3,1)$  & 2 & 2 \\
\hline
$\Gamma_{5}$ & $(0,1,1,1,2,1)$  & 2 & 2 \\
\hline
$\Gamma_{6}$ & $(0,0,2,2,1,1)$  & 3 & 1 \\
\hline
$\Gamma_{7}$ & $(4,2,1,3,2,1)$  & 3 & 4 \\
\hline
$\Gamma_{8}$ & $(0,0,1,1,2,1)$  & 3 & 2 \\
\hline
$\Gamma_{9}$ & $(1,0,1,1,1,0)$  & 3 & 2 \\
\hline
$\Gamma_{10}$ & $(2,0,1,1,1,0)$  & 4 & 2 \\
\hline
$\Gamma_{11}$ & $(0,0,1,1,1,0)$  & 4 & 2 \\
\hline
$\Gamma_{12}$ & $(2,2,0,2,1,1)$  & 4 & 2 \\
\hline
$\Gamma_{13}$ & $(0,0,0,0,1,1)$  & 4 & 1 \\
\hline
\end{tabular}
\end{center}

\noindent Further, we will calculate explicitly the singular points of complex hypersurfaces $\tilde{s}(t)_{\Gamma}$ for every face $\Gamma \subset \Delta$ from the table above. Also we will consider only one face from every orbit because argumentation for other faces is always similar.

At first, note that the face $\Gamma_{1}$ is a parallelogram and corresponding truncated scalar curvature has coefficients $-b_2,-b_3,-b_3,-b_4$, so according to Lemma \ref{prlm} the corresponding discriminant condition can be written as the equality $b_2 b_5-b_3^2=0$ which is equivalent to the condition (\ref{eqth2}).

It is easy to see that the systems corresponding to faces $\Gamma_3, \Gamma_4, \Gamma_5, \Gamma_6, \Gamma_8$ are inconsistent due to that $a_1+2b_1>0$ and $a_2 + 2 b_2 > 0$.

Consider the face $\Gamma_2$ and the following change of variables $x=t_3 t_4^{-1}, \ y=t_1 t_2^{-1}$. The system of equations takes the following form
\begin{numcases}{\Gamma_{2}:}
\label{2.1} t_1 \tilde{s}_{\Gamma_{2}}(t)=a_1+a_2 y - (b_1+b_2 y)\left(x + x^{-1} \right)=0, \\
\label{2.2} (b_1+b_2 y)\left(1 - x^{-2} \right)=0, \\
\label{2.3} a_2 - b_2 \left(x+ x^{-1} \right) = 0,
\end{numcases}
\noindent From (\ref{2.3}) we immediately obtain that $x=\frac{a_2 \pm \sqrt{a_2^{2}-4 b_2^2}}{2 b_2}$. Using the inequality $a_2 \pm 2 b_2 > 0$ one can show that $x \neq \pm 1$, so $y=-\frac{b_1}{b_2}$. Substituting the expressions for $x$ and $y$ to the equation (\ref{2.1}) we find out the consistency condition $a_1 b_2 - a_2 b_1 = 0$. This condition is equivalent to (\ref{eqth1}).

Consider the face $\Gamma_7$ and the following change of variables $x= t_4^{2} t_1^{-1} t_2^{-1}, \ y= t_4 t_5 t_1^{-1} t_6^{-1}, \ z= t_3 t_4 t_5^{-2}$. The system of equations takes the following form
\begin{numcases}{\Gamma_{7}:}
\label{7.1} -t_1^{-1} t_3 t_4 \tilde{s}_{\Gamma_{7}}(t) = b_1 + b_2 x + b_3 y (1+z) + b_4 z + b_5 y^2 z x^{-1} = 0, \\
\label{7.2} b_2 - b_5 y^2 z x^{-2} = 0, \\
\label{7.3} b_3 (1+z) + 2 b_5 yz x^{-1} = 0, \\
\label{7.4} b_3 y + b_4 + b_5 y^2 x^{-1} = 0,
\end{numcases}
\noindent From the equation (\ref{7.2}) we have $y^2 z x^{-2}=\frac{b_2}{b_5}$. On the other hand, from (\ref{7.3}) it follows that $yz x^{-1}=-\frac{b_3}{2 b_5}(1+z)$, so we obtain the equality $2 b_2 x = - b_3 (1+z) y$. Keeping in mind the derived equalities and (\ref{7.1}) we find out that $z=-\frac{b_1}{b_4}$. From this one can deduce that $y x^{-1} = \frac{b_3}{2 b_5}(\frac{b_4}{b_5}-1)$ and $y^2 x^{-2}=-\frac{b_2 b_4}{b_1 b_5}$ simultaneously. Since $n_1,~n_2,~n_3 \in \mathbb{N}$ we obtain that $y x^{-1} \in \mathbb{R}$ and $(y x^{-1})^{2} < 0$ which is impossible. Therefore, the system for the face $\Gamma_{7}$ is inconsistent.

Consider the face $\Gamma_9$ and the following change of variables $x= t_2 t_6^{-1},~ y= t_3 t_4^{-1},~ z= t_3 t_5^{-1}$. The system of equations takes the following form
\begin{numcases}{\Gamma_{9}:}
\label{9.1} t_2 \tilde{s}_{\Gamma_{9}}(t)=a_2 + a_5 x - b_2 \left( y + y^{-1} \right) - b_3 x \left( z + z^{-1} + z y^{-1} + y z^{-1} \right) - b_5 x^2 = 0, \\
\label{9.2} a_5 - b_3 \left( z + z^{-1} + z y^{-1} + y z^{-1} \right) - 2 b_5 x = 0, \\
\label{9.3} b_2 \left( 1 - y^{-2} \right) + b_3 x \left( z^{-1} - z y^{-2} \right) = 0, \\
\label{9.4} b_3 x \left( 1 - z^{-2} + y^{-1} - y z^{-2} \right) = 0.
\end{numcases}
\noindent In case of $y=\pm 1$, from  linear combination of the equations (\ref{9.1}) and (\ref{9.2}) it follows that $x^2=-\frac{a_2 \mp b_2}{b_5} < 0$. On the other hand, from (\ref{9.3}) we have $z=\pm 1$, so $x \in \mathbb{R}$ due to (\ref{9.2}), but it is impossible. Consequently, $y \neq \pm 1$. It means that $y=z^2$ and the equation (\ref{9.3}) gives us the equality $b_2 (z+\frac{1}{z})+b_3 x =0$. Then from (\ref{9.2}) we have $x=-\frac{a_5 b_2}{2(b_3^2- b_2 b_5 )}$ ($b_3^2-b_2 b_5=0$ leads to the contradiction $a_5=0$). From the previous equality we obtain $z+\frac{1}{z}=\frac{a_5 b_3}{2(b_3^2-b_2 b_5)}$. From this it immediately follows that $y+\frac{1}{y}=z^2 + \frac{1}{z^2}=(z+\frac{1}{z})^2-2=\frac{a_5^2 b_3^2}{4(b_3^2-b_2 b_5)^2} - 2$. After substitution of these expressions into linear combination of the equations (\ref{9.1}) and (\ref{9.2}), namely, $a_2 - b_2 (y+\frac{1}{y})+ b_5 x^2=0$ we have the condition of consistency $4(b_3^2-b_2 b_5)(a_2+2b_2)=a_5^2 b_2$ which is equivalent to
\begin{equation}\label{eqth4}
\Big( 2 n_1 ( 2 n_3 + 1) - (n_2 - 1)^2 \Big)(2 n_1 + n_2 - 1) - \Big( 2(n_1+n_2+n_3) - 1 \Big)^2 (2 n_3 + 1) = 0.
\end{equation}

\noindent After easy manipulations we get
\begin{equation}
(2 n_3 + 1) \Big( 2 n_1 (2n_1 + n_2 - 1) - (2 (n_1 + n_2 + n_3) - 1)^2 \Big) - (n-1)^2 ( 2 n_1 + n_2 -1) = 0.
\end{equation}

\noindent On the other hand, it is easy to see that
\begin{equation*}
2 n_1 (2n_1 + n_2 - 1) - (2 (n_1 + n_2 + n_3) - 1)^2 =
\end{equation*}
\begin{equation}
=- 2 n_1 (3 n_2 - 1) - 4 n_2 (n_2 - 1) - 4 n_3 (n_3 - 1) - 8 (n_1 + n_2) n_3 - 1 < 0.
\end{equation}

\noindent Therefore, the equality (\ref{eqth4}) cannot be true, so the corresponding system is inconsistent.

Consider the face $\Gamma_{10}$ and the following change of variables $x= t_6 t_2^{-1},~ y= t_5 t_3^{-1},~ z= t_1 t_6 t_3^{-1} t_4^{-1},~ w= t_3 t_4 t_5^{-2}$. The system of equations takes the following form
\begin{numcases}{\Gamma_{10}:}
t_6 \tilde{s}_{\Gamma_{10}}(t)=a_5 + a_2 x - b_1 z - b_2 x \left( y^2 w + y^{-2} w^{-1} \right) - \nonumber \\
\label{10.1} - b_3 \left( y + y^{-1} + y w + y w^{-1} \right) - b_4 zw - b_5 x^{-1} = 0, \\
\label{10.2} a_2 - b_2 \left( y^2 w + y^{-2} w^{-1} \right) + b_5 x^{-2} = 0, \\
\label{10.3} b_2 x \left( 2 y w - 2 y^{3} w^{-1} \right) + b_3 \left( 1 - y^{-2} + w - y^{-2} w^{-1} \right) = 0, \\
\label{10.4} b_1 + b_4 w = 0, \\
\label{10.5} b_2 x \left( y^2  - y^{-2} w^{-2} \right) + b_3 \left( y - y^{-1} w^{-2} \right) + b_4 z = 0,
\end{numcases}
\noindent From the equation (\ref{10.4}) we immediately obtain the equality $w = - \frac{b_1}{b_4} < 0$. Let us perform the following transformations of the equation (\ref{10.3})
$$
b_2 x \left( 2 y w - 2 y^{-3} w^{-1} \right) + b_3 \left( 1 - y^{-2} + w - y^{-2} w^{-1} \right) = 2 b_2 x y w \left(1 - y^{-2} w^{-1} \right) \left( 1 + y^{-2} w^{-1} \right) +$$
$$+ b_3 (w+1) \left( 1 - y^{-2} w^{-1} \right)=\left( 1 - y^{-2} w^{-1} \right) \left( 2 b_2 x y w \left( 1 + y^{-2} w^{-1} \right) + b_3 (w+1) \right)=0.
$$
\noindent The case $y^2 w = 1$ is impossible because then $y^2 < 0$ and, therefore, $\text{Re}~ y=0$. Also $x^2=-\frac{b_5}{a_2 -2 b_2} < 0$ and, correspondingly, $\text{Re}~x=0$. Then after taking real part of linear combination of the equations (\ref{10.1}) and (\ref{10.4}) we obtain contradiction $a_5=0$. Also the case $y^2 w = -1$ is impossible because then $w=-1$ and $y^2=1$, so $z=0$ due to the last equation. Therefore, $y^2 w \neq \pm 1$ and $w \neq -1$. Then we have the following equality
$$
2 b_2 x y w \left( 1 + y^{-2} w^{-1} \right) + b_3 (w+1) = 0.
$$
\noindent From this we obtain
$$
x = - \frac{b_3}{2 b_2} \frac{y(w+1)}{y^2 w + 1} \Rightarrow x^{-2} = \frac{4 b_2^2}{b_3^2} \frac{(y^2 w + 1)^2}{y^2 (w+1)^2} = \frac{8 b_2^2}{b_3^2} \frac{w}{(w+1)^2} + \frac{4 b_2^2}{b_3^2} \frac{w}{(w+1)^2} \left( y^2 w + y^{-2} w^{-1} \right).
$$
\noindent On the other hand, from (\ref{10.2}) we have the equality $b_2 \left( y^2 w + y^{-2} w^{-1} \right) = b_5 x^{-2} + a_2$. After substitution of this equality into the previous one we obtain
$$
\left( \frac{b_3^2}{4 b_2} (w+1)^2 w^{-1} - b_5 \right) x^{-2} = a_2 + 2 b_2 > 0.
$$
\noindent From this it follows that $x^2 < 0$, so $\text{Re}~x=0$. Besides, by using the equation (\ref{10.2}) we obtain $y^2 w + y^{-2} w^{-1} \in \mathbb{R}$. Further from linear combination of the equations (\ref{10.3}) and (\ref{10.5}) it follows that
$$
b_3 \left( y w^{-1} - y^{-1} w^{-1} - y + y^{-1} w^{-2} \right) - 2 b_4 z = 0,
$$
\noindent from which we immediately have the following
$$
z=\frac{b_3}{2 b_4} (1-w) w^{-1} \left( y + y^{-1} w^{-1} \right) \Rightarrow z^2=\frac{b_3^2}{4 b_4^2} (1-w)^2 w^{-3} \left( y^2 w + y^{-2} w^{-1} + 2 \right) \in \mathbb{R}.
$$
\noindent Now consider the case $y^2 w + y^{-2} w^{-1} < -2$. Then it is true that $y^2 w < 0$ and $y^2>0$ because the equation $t + t^{-1}=C,~ C<-2$ always has two negative roots, so $\text{Im}~y=0$. After taking imaginary part of linear combination of the equations (\ref{10.1}) and (\ref{10.4}) we obtain the following equality
$$
a_2 x - b_2 x \left( y^2 w + y^{-2} w^{-1} \right) - b_5 x^{-1} = 0,
$$
\noindent linear combination of which with (\ref{10.2}) gives contradiction $b_5 x^{-1}=0$. The case $y^2 w + y^{-2} w^{-1} = -2$ is also impossible because then $z=0$. At last, if $y^2 w + y^{-2} w^{-1} > -2$ then $\left( y + y^{-1} + y w + y^{-1} w^{-1} \right)^2 = (1+w)^2 (y+ y^{-1} w^{-1} )^2= (1+w)^2 w^{-1} (y^2 w + y^{-2} w^{-1} + 2) < 0$. After taking real part of linear combination of the equations (\ref{10.1}) and (\ref{10.2}) we obtain contradiction $a_5=0$. Consequently, the system for the face $\Gamma_{10}$ is inconsistent.

Consider the face $\Gamma_{11}$ and the following change of variables $x= t_6 t_1^{-1},~ y= t_6 t_2^{-1},~ z= t_4 t_5^{-1},~ w= t_3 t_4^{-1}$. The system of equations takes the following form
\begin{numcases}{\Gamma_{11}:}
t_6 \tilde{s}_{\Gamma_{12}}(t)= a_1 x + a_2 y + a_5  - (b_1 x + b_2 y) \left( w + w^{-1} \right) - \nonumber \\
\label{11.1} - b_3 \left( z + z^{-1} + zw + z^{-1} w^{-1} \right) - b_5 y^{-1} = 0, \\
\label{11.2} a_1 - b_1 \left( w + w^{-1} \right) = 0, \\
\label{11.3} a_2 - b_2 \left( w + w^{-1} \right) + b_5 y^{-2} = 0, \\
\label{11.4} b_3 (1 - z^{-2} + w - z^{-2} w^{-1} ) = 0, \\
\label{11.5} (b_1 x + b_2 y) \left( 1 - w^{-2} \right) + b_3 \left( z - z^{-1} w^{-2} \right) = 0,
\end{numcases}
\noindent From the 2nd equation we obtain that $w + w^{-1} = \frac{a_1}{b_1} > 2$ and, consequently, $w \in \mathbb{R},~ w>0$. After substitution of this equality to the 3rd equation we have that $y^{-2}=\frac{a_1 b_2 - a_2 b_1}{b_1 b_5} \in \mathbb{R}$. Since $1- z^{-2}+w- z^{-2} w^{-1} = (1+w)(1 - z^{-2} w^{-1} )$ then from the 4th equation we obtain that $z^2=w^{-1}>0$. Then $z+ z^{-1} + zw + z^{-1} w^{-1}=2(z + z^{-1} )$ and from the linear combination of the 1st, 2nd and 3rd equation we obtain the following equality
$$
a_5 - 2 b_3 \left( z + z^{-1} \right) - 2 b_5 y^{-1} = 0,
$$
\noindent which after rising to square takes the form
$$
y^{-1}=\frac{1}{4 a_5 b_1 b_5} \left( a_5^2 b_1 + 4 b_5 (a_1 b_2 - a_2 b_1) - 4 b_3^2 (a_1 + 2 b_1) \right).
$$
\noindent Whence we immediately obtain the following condition of inconsistency of the system
\begin{equation}\label{eq7}
16 a_5^2 b_1 b_5 (a_1 b_2 - a_2 b_1) = (a_5^2 b_1 + 4 b_5 (a_1 b_2 -a_2 b_1) - 4 b_3^2 (a_1 + 2 b_1))^2.
\end{equation}
\noindent that can be written in the form of (\ref{eqth7}). When it holds we have two metrics.

Consider the face $\Gamma_{12}$ and the following change of variables $x= t_4 t_1^{-1},~ y=\ t_1 t_6 t_4^{-1} t_5^{-1},~ z= t_2 t_5 t_4^{-1} t_6^{-1},~ w= t_3 t_4 t_5^{-1} t_6^{-1}$. The system of equations takes the following form
\begin{numcases}{\Gamma_{12}:}
t_3 \tilde{s}_{\Gamma_{11}}(t)= a_3 - b_1 \left( x + x^{-1} \right) - b_2 \left( xyz + x^{-1} y^{-1} z^{-1} \right) - \nonumber \\
\label{12.1} - b_3 \left( xy + x^{-1} y^{-1} + w \right) - b_4 yw - b_5 zw = 0, \\
\label{12.2} b_1 \left( x - x^{-1} \right) + b_2 \left( xyz - x^{-1} y^{-1} z^{-1} \right) + b_3 \left( xy - x^{-1} y^{-1} \right) = 0, \\
\label{12.3} b_2 \left( xyz - x^{-1} y^{-1} z^{-1} \right) + b_3 \left( xy - x^{-1} y^{-1} \right) + b_4 yw = 0, \\
\label{12.4} b_2 \left( xyz - x^{-1} y^{-1} z^{-1} \right) + b_5 zw =0, \\
\label{12.5} b_3 + b_4 y + b_5 z =0,
\end{numcases}
\noindent From (\ref{12.5}) we immediately obtain $y=-\frac{b_3+b_5 z}{b_4}$. Consider the equality obtained as the sum of the equations (\ref{12.2}), (\ref{12.3}) and (\ref{12.4}) with coefficients $b_5 z$, $(-1) b_5 z$ and $b_4 y$, correspondingly. This equality has the following form
$$
b_1 b_5 z \left( x- x^{-1} \right) + b_2 b_4 y \left( xyz - x^{-1} y^{-1} z^{-1} \right) = 0.
$$
Multiplying this equality by $zx$ we obtain
$$
b_1 b_5 z^2 \left( 1 - x^{-2} \right) + b_2 b_4 \left( y^2 z^2 - x^{-2} \right) = 0,
$$
\noindent from which we have the following expression
$$
x^{-2} = \frac{(b_1 b_5 + b_2 b_4 y^2) z^2}{b_1 b_5 z^2 + b_2 b_4}.
$$
\noindent Moreover, after multiplication of the equation (\ref{12.2}) by $\frac{yz}{x}$ we obtain the following
$$
b_1 yz \left( 1 - x^{-2} \right) + b_2 \left( y^2 z^2 - x^{-2} \right) + b_3 z \left( y^2 - x^{-2} \right) = 0.
$$
\noindent After substitution the expression for $x^{-2}$ we obtain the equality of the following form
$$
b_1 yz \left( 1 - \frac{(b_1 b_5 + b_2 b_4 y^2) z^2}{b_1 b_5 z^2 + b_2 b_4} \right) + b_2 \left( y^2 z^2 - \frac{(b_1 b_5 + b_2 b_4 y^2) z^2}{b_1 b_5 z^2 + b_2 b_4} \right) + b_3 z \left( y^2 - \frac{(b_1 b_5 + b_2 b_4 y^2) z^2}{b_1 b_5 z^2 + b_2 b_4} \right) = 0,
$$
\noindent which after substitution $y=-\frac{b_3+b_5 z}{b_4}$ and additional transformations takes the form
$$
p_1(z)=2 b_1 b_2 b_5^3 z^5 + b_3 b_5^2 (5 b_1 b_2 + b_1 b_5 - b_2 b_4 ) z^4 + b_3^2 b_5 (4 b_1 b_2 + 2 b_1 b_5 - 2 b_2 b_4 ) z^3 +$$ $$+ b_3 (b_1 b_2 b_3^2 + b_1 b_3^2 b_5 - b_1 b_4^2 b_5 - b_2 b_3^2 b_4 + b_2 b_4 b_5^2) z^2 + 2 b_2 b_4 b_5 ( b_3^2 - b_1 b_4 ) z + b_2 b_3 b_4 ( b_3^2 - b_1 b_4) =0.
$$
\noindent On the other hand, linear combination of the equations (\ref{12.1}), (\ref{12.2}) and (\ref{12.5}) gives us the following relation
$$
a_3 - 2 b_1 x - 2 b_2 xyz - b_3 xy = 0 \Rightarrow a_3 x^{-1} = 2 b_1 + 2 b_2 yz + 2 b_3 y,
$$
\noindent which squaring of both sides and multiplication on the denominator of the left hand side we obtain the following
$$
a_3^2 (b_1 b_5 + b_2 b_4 y^2) z^2 = 4 (b_1 + b_2 yz + b_3 y)^2 (b_1 b_5 z^2 + b_2 b_4).
$$
\noindent After substitution into this equality of the expression for $y$ and additional transformations it has the form
$$
p_2 (z) = 4 b_1 b_2^2 b_5^3 z^6 + 8 b_1 b_2 b_3 b_5^2 ( b_2 + b_5) z^5 +
$$
$$
+b_5 (4 b_2^3 b_4 b_5 + 4 b_1 b_3^2 b_5^2 + 16 b_1 b_2 b_3^2 b_5 + 4 b_1 b_2^2 b_3^2  - 8 b_1^2 b_2 b_4 b_5 - a_3^2 b_2 b_4 b_5 ) z^4 + $$ $$+ 2 b_3 b_5 (4 b_1 (b_2 + b_5) (b_3^2 - b_1 b_4 )+ b_2 b_4 ( 4 b_2 (b_2 + b_5) - a_3^2 )) z^3 + $$ $$+ (16 b_2^2 b_3^2 b_4 b_5 + 4 b_2 b_3^2 b_4 b_5^2 + 4 b_2^3 b_3^2 b_4 + 4 b_1 b_3^4 b_5 + 4 b_1^3 b_4^2 b_5 - 8 b_1^2 b_3^2 b_4 b_5  - 8 b_1 b_2^2 b_4^2 b_5 - a_3^2 b_1 b_4^2 b_5 - a_3^2 b_2 b_3^2 b_4  ) z^2 +$$ $$+ 8 b_2 b_3 b_4 (b_3^2 - b_1 b_4) (b_2 + b_5) z + 4 b_2 b_4 (b_3^2 - b_1 b_4)^2 = 0.
$$
\noindent Therefore, the condition of consistency has the following form
\begin{equation}\label{eq6}
\text{Res}(p_1,p_2)=0,
\end{equation}
\noindent where the left-hand side is the resultant of polynomials $p_1$ and $p_2$. 

Consider the face $\Gamma_{13}$ and the following change of variables $x= t_6 t_5^{-1},~ y= t_3 t_4^{-1},~ z= t_4 t_1^{-1},~ w= t_4 t_2^{-1}$. The system of equations takes the following form
\begin{numcases}{\Gamma_{13}:}
t_3 \tilde{s}_{\Gamma_{13}}(t)= a_3 + (a_3 + a_1 z + a_2 w) y - b_1 \left( z + z^{-1} + y^2 z \right) - \nonumber \\
\label{13.1} - b_2 \left( w + w^{-1} + y^2 w \right) - b_3 (1+y) \left( x + x^{-1} \right) = 0, \\
\label{13.2} b_3 (1+y) \left( 1 - x^{-2} \right) = 0, \\
\label{13.3} a_3 + a_1 z + a_2 w - 2 b_1 y z - 2 b_2 y w - b_3 \left( x + x^{-1} \right) = 0, \\
\label{13.4} a_1 y - b_1 \left( 1 - z^{-2} + y^2 \right) = 0, \\
\label{13.5} a_2 y - b_2 \left( 1 - w^{-2} + y^2 \right) = 0,
\end{numcases}
\noindent From the equation (\ref{13.2}) we immediately obtain either $y=-1$ or $x=\pm1$. Consider the case when $x = \pm 1$. Then from the equations (\ref{13.4}) and (\ref{13.5}) we obtain the following expressions
$$
b_1 z^{-2} = b_1 (y^2 + 1) - a_1 y, \ \ \ b_2 w^{-2}= b_2 (y^2 + 1) - a_2 y.
$$
\noindent After squaring the equation (\ref{13.3}) twice and corresponding substitutions this equation takes the form
$$
q_{1}^{\mp}(y)=(a_3 \mp 2 b_3)^4 (b_1 (y^2 + 1) - a_1 y)^2 (b_2 (y^2 + 1) - a_2 y)^2 + b_1^2 (2 b_1 y - a_1)^4 (b_2 (y^2 + 1) - a_2 y)^2 + $$ $$+ b_2^2 (2 b_2 y - a_2)^4 (b_1 (y^2 + 1) - a_1 y)^2 - 2 b_1 (a_3 \mp 2 b_3)^2 (2 b_1 y - a_1)^2 (b_1 (y^2 + 1) - a_1 y) (b_2 (y^2 + 1) - a_2 y)^2 - $$ $$- 2 b_2 (a_3 \mp 2 b_3)^2 ( 2 b_2 y - a_2)^2 (b_1 (y^2 + 1) - a_1 y)^2 (b_2 (y^2 + 1) - a_2 y) -$$ $$- 2 b_1 b_2 (2 b_1 y - a_1)^2 (2 b_2 y - a_2)^2 (b_1 (y^2 + 1) - a_1 y) (b_2 (y^2 + 1) - a_2 y) = 0.
$$
\noindent On the other hand, from linear combination of the equations (\ref{13.1}), (\ref{13.2}), (\ref{13.4}) and (\ref{13.5}) we obtain the following equality
$$
a_3 ( 1 + y) = 2 \left( b_1 z^{-1} + b_2 w^{-1} \right),
$$
\noindent which after squaring both sides twice and corresponding substitutions takes the form
$$
q_2(y)= a_3^4 (1+y)^4 - 8 a_3^2 (1+y)^2 \left( (b_1^2 + b_2^2) (y^2 + 1) -  (a_1 b_1 + a_2 b_2 ) y \right) + $$ $$+ 16 \left( (b_1^2 - b_2^2) (y^2 + 1) + (a_2 b_2 - a_1 b_1) y \right)^2 = 0.
$$
\noindent Therefore, we have the following two conditions of consistency
\begin{equation}\label{eq9}
\text{Res}(q_{1}^{+}, q_2)=0, \ \ \ \ \ \text{Res}(q_{1}^{-}, q_2)=0.
\end{equation}

Now let $y=-1$. Then from (\ref{13.4}), (\ref{13.5}) we have that $z^{-2}=2+\frac{a_1}{b_1}, \ w^{-2}=2+\frac{a_2}{b_2}$. After substitution the values of variables $y,z,w$ into (\ref{13.1}) we obtain the following equality for parameters
\begin{equation*}
a_1 z + a_2 w + b_1 z \left( 2 + z^{-2} \right) + b_2 w \left( 2 + w^{-2} \right) = 0
\end{equation*}
\begin{equation*}
a_1 z + a_2 w + b_1 z \left( 4 + \frac{a_1}{b_1} \right) + b_2 w \left( 4 + \frac{a_2}{b_2} \right) = 0
\end{equation*}
\begin{equation*}
( a_1 + 2 b_1 ) z + ( a_2 + 2 b_2 ) w = 0
\end{equation*}
\begin{equation*}
\frac{( a_1 + 2 b_1 )^2}{2+\frac{a_1}{b_1}}  =  \frac{( a_2 + 2 b_2 )^2}{2+\frac{a_2}{b_2}}
\end{equation*}
\begin{equation}\label{eq13}
b_1 ( a_1 + 2 b_1 ) =  b_2 ( a_2 + 2 b_2 )
\end{equation}
\noindent It can be shown that the last equality (\ref{eq13}) is equivalent to the condition $n_1- n_2 = 0$.
\hfill$\Box$

\section{Inonu-Wigner contractions of flag manifolds}

Now describe construction of Inonu-Wigner contraction of a flag manifold $G/H$. Let $\Delta$ be corresponding Newton polytope and $\Gamma \subset \Delta$ is some proper face of it. As it was already mentioned, the face $\Gamma$ can be defined by the normal vector $f \in \mathbb{R}^{n}$ with integer coefficients. In its turn, this vector $f$ defines the filtration of the Lie algebra $\mathfrak{g}$ by the following way
\begin{equation}
F_{a} \mathfrak{g} = \mathfrak{h} \oplus \left( \bigoplus\limits_{f(j) \leq a} \mathfrak{m}_{j} \right), \ \ \ a \in \mathbb{Z}.
\end{equation}

\noindent As usual we can now take the graded Lie algebra $\mathfrak{g}_{f}$ corresponding to this filtration as follows
\begin{equation}
\mathfrak{g}_{\Gamma}= \bigoplus\limits_{a \in \mathbb{Z}} \mathfrak{g}^{(a)}_{\Gamma}, \ \ \text{where} \ \ \mathfrak{g}^{(a)}_{\Gamma}=F_{a}\mathfrak{g} / \bigcup\limits_{b<a} F_{b} \mathfrak{g}.
\end{equation}

\noindent Note that we have the inclusion $\mathfrak{h} \subset \mathfrak{g}_{\Gamma}$, so there exists the Lie group $G_{\Gamma}$ with Lie algebra $\mathfrak{g}_{\Gamma}$ and the subgroup $H \subset G_{\Gamma}$. Therefore, we can consider the homogeneous manifold $M_{\Gamma}=G_{\Gamma}/H$ which is called Inonu-Wigner contraction of manifold $M$ by the face $\Gamma \subset \Delta$.

It is important that Newton polytope of the contraction $M_{\Gamma}$ is very face $\Gamma$ and, moreover, its scalar curvature is truncation of the scalar curvature of $M$ on the face $\Gamma$. The same can be said also for corresponding Einstein equations. Also it is easy to see that Einstein metrics on $M_{\Gamma}$ have vanishing scalar curvature and all its derivatives, so the Ricci tensor vanish. Consequently, Einstein metrics on the contraction $M_{\Gamma}$ are Ricci-flat, so we have the following theorem.

\begin{Teo}[{\cite[Th. 3.1]{GrUMN}}]
\textit{Flag manifold $M$ has positive defect $\nu(M) - \mathcal{E}(M) > 0$ if and only if there exists invariant complex Einstein metric $g^{M_{\Gamma}}$ on the contraction $M_{\Gamma}$ of manifold $M$ by some proper face $\Gamma \subset \Delta(M)$. Moreover, such Einstein metrics $g^{M_{\Gamma}}$ are Ricci-flat.}
\end{Teo}

It is worth to note that if the face $\Gamma$ is a parallelogram and there exist Einstein metrics on the contraction $M_{\Gamma}$, then one can choose real metrics from these Einstein metrics. Moreover, in most cases one can choose Lorentzian metrics. For example, we have the following theorem for the flag manifold $M_{n_1, n_2, n_3}$ considered in the previous section, where the faces $\Gamma_{1}^{i,j}$ are elements of the orbit of the face $\Gamma_{1}$ which can be defined by the normal vectors $f_{1}^{1,1}=(4,2,2,4,3,1), f_{1}^{1,2}=(4,2,4,2,3,1), f_{1}^{2,1}=(2,4,2,4,1,3), f_{1}^{2,2}=(2,4,4,2,1,3)$.

\begin{Teo}
\textit{If the equality $8 n_i (2 n_3 + 1) - (n_j - 1)^2 = 0$ is true, where $(i,j)=(1,2)$ or $(2,1)$, then contractions of the manifold $M_{n_1,n_2,n_3}$ by the faces $\Gamma_{1}^{i,1}$ and $\Gamma_{1}^{i,2}$ with the real metrics given by
$$g^{M_{\Gamma_{1}^{1,1}}}=\left( t_1, t_2, t_3,  \frac{(n_2 - 1)^2}{4 n_1^2} t_2^2 t_3, \frac{n_2 - 1}{2 n_1} t_2 t_3, -1 \right), \ \ t_1, t_2, t_3 \in \mathbb{R}_{>0},$$
$$g^{M_{\Gamma_{1}^{1,2}}}=\left( t_1, t_2, t_3,  \frac{4 n_1^2}{(n_2 - 1)^2} t_2^{-2} t_3, \frac{2 n_1}{n_2 - 1} t_2^{-1} t_3, -1 \right), \ \ t_1, t_2, t_3 \in \mathbb{R}_{>0},$$
$$g^{M_{\Gamma_{1}^{2,1}}}=\left( t_1, t_2, t_3, t_3^{-1}, \frac{n_1 - 1}{2 n_2} t_1 t_3, -1 \right), \ \ t_1, t_2, t_3 \in \mathbb{R}_{>0},$$
$$g^{M_{\Gamma_{1}^{1,1}}}=\left( t_1, t_2, t_3, t_3^{-1}, \frac{n_1 - 1}{2 n_2} t_2 t_3^{-1}, -1 \right), \ \ t_1, t_2, t_3 \in \mathbb{R}_{>0},$$
are Ricci-flat manifolds.}
\end{Teo}

\noindent \textit{Proof:} As in the proof of Theorem \ref{number}, consider the discriminant condition for the face $\Gamma_1^{1,1}$ in the following non-homogeneous form
\begin{numcases}{\Gamma_{1}^{1,1}:}
\label{1.1} -t_2 t_3 t_4^{-1} \tilde{s}_{\Gamma_{1}}(t)=b_2+b_3(x+y)+b_5 x y=0, \\
\label{1.2} b_3 + b_5 y = 0, \\
\label{1.3} b_3 + b_5 x = 0.
\end{numcases}
\noindent where $x= t_2 t_3 t_5^{-1} t_6^{-1},~ y=t_2 t_5 t_4^{-1} t_6^{-1}$. From the equations (\ref{1.2}) and (\ref{1.3})  we find out that $x=-\frac{b_3}{b_5}, y=-\frac{b_3}{b_5}$. So if the condition of consistency $b_2 b_5 - b_3^2=0$ which is equivalent to the equality mentioned in the statement of the theorem for $(i,j)=(1,2)$ holds, then the contraction $M_{\Gamma_1^{1,1}}$ with the metric $g^{M_{\Gamma_1^{1,1}}}=(t_1,..,t_6)$ such that $t_2 t_3 t_5^{-1} t_6^{-1} = t_2 t_5 t_4^{-1} t_6^{-1}= -\frac{b_3}{b_5} = - \frac{2 n_1}{n_2 - 1}$ is Ricci-flat. Moreover, we can fix $f_6=-1$ and variables $t_4, t_5$ express through $t_1, t_2, t_3$. Therefore, we obtain family of real Ricci-flat metrics mentioned in the statement of the theorem parametrized by $t_1, t_2, t_3 \in \mathbb{R}_{>0}$. For other faces argumentation is similar.
\hfill$\Box$

\noindent \textbf{Remark}: In fact, these homogeneous manifolds are diffeomorphic to Euclidean spaces.

\end{document}